%% file: rationalTrig.tex
\documentclass[11pt]{article}
\usepackage{geometry}                % See geometry.pdf to learn the layout options. There are lots.
\usepackage{amssymb}
\usepackage{epstopdf}
\usepackage{xr}		% external references
\input{GunnMac}

\date{January 13, 2014\footnote{slightly updated since then}}

\renewcommand{\vec}[1]{\mathbf{#1}}

\newcommand{\gTh}{\cite{gunnThesis}\xspace}

\newcommand{\tarea}{\ensuremath{\mathcal{A}}}
\begin{document}
\title{Rational trigonometry via projective geometric algebra}
\author{Charles G. Gunn, Ph. D.\footnote{projgeom@gmail.com, Raum+Gegenraum, Falkensee, Germany}
}

\maketitle
%\label{ch:conformal}
%\section{Comparison with conformal model of euclidean geometry}

\begin{abstract}
%We show that main results of \cite{wildberger2005} regarding rational trigonometry and  universal geometry can be expressed in  projective geometric algebra as developed in \gTh and \cite{gunn2011}.    This provides a coordinate-free representation for these results with promising compactness and fidelity.  
We show that main results of rational trigonometry and universal geometry (\cite{wildberger2005}) can be succinctly expressed using projective geometric algebra (PGA) (\cite{gunn2011}, \cite{gunn2017a}).   In fact, the PGA representation exhibits distinct advantages over the original vector-based approach.   These include the advantages intrinsic to geometric algebra: it is coordinate-free, treats lines and points in a unified framework, and handles many special cases in a uniform and seamless fashion.  It also reveals structural patterns not visible in the original formulation, for example, the exact duality of spread and quadrance.  The current article handles only a representative (euclidean) subset of the full content of Wildberger's work, but enough we believe to establish the value of this approach for further development, also for the noneuclidean cases.
\end{abstract}

\section{Introduction}

We apply the geometric algebra $\pdclal{2}{0}{1}$ (\cite{gunn2017b}) to the geometric entities of \cite{wildberger2005}.  This algebra is the 2D euclidean variant of \emph{projective geometric algebra}. In the sequel we sometimes shorten this to {PGA}.  This article is restricted to the euclidean case, but PGA is metric-neutral and the treatment presented here can be easily extended to noneuclidean metrics.  

\subsection{Crash course in 2D PGA}\label{sec:crashcourse}
We begin with an unsystematic listing of some key aspects of $\pdclal{2}{0}{1}$ which are needed for what follows. Consult \cite{gunn2017b} for details. 
\begin{itemize}
\item The geometric, or Clifford algebra, is built atop the Grassmann, or exterior algebra, and inherits its graded structure.  The Grassmann algebra has an anti-symmetric product, the wedge product, which is additive on grade: the product of a $k$-vector and an $l$-vector is a $k+l$-vector (when it isn't zero), which geometrically corresponds to the \emph{join} of two elements.  There is also a dual Grassmann algebra in which 1-vectors represent dual vectors (in dimension 2, lines); in this algebra, the wedge product represents the \emph{meet} operator. To differentiate these two wedge products, we use $\wedge$ to represent meet and $\vee$ to represent join.  See the references to understand how both these operations can be accessed from within a single Grassmann algebra.  
\item  For a multivector $\vec{X}$, $\langle \vec{X} \rangle_k$ represents the grade-$k$ part of $\vec{X}$, so $\vec{X} = \sum_{k=0}^3\langle \vec{X} \rangle_k$.  
\item The basis of our Clifford algebra is the dual, projectivized Grassmann algebra:
\begin{itemize}
\item A line $ax+by+c=0$ maps to the 1-vector $c\e{0}+a\e{1}+b\e{2}$.
\item The point $(x,y)$  maps to the 2-vector $\EE{0} +x\EE{1}+y\EE{2}$, where the basis 2-vectors are defined as: $\EE{i} = \e{j} \e{k}$ for $ijk$ a cyclic permutation of $012$.  
\item The point $(x,y,z)$ (in homogeneous coordinates) maps to the 2-vector $z\EE{0} +x\EE{1}+y\EE{2}$.   A point with $z\neq0$ can be normalized to have $z=1$.  
\end{itemize}
\item One obtains a geometric algebra by attaching an inner product (a bilinear symmetric form)  to the underlying vector space.  The geometric product is then defined for 1-vectors as
\[ \vec{l}\vec{m} = \vec{l}\cdot \vec{m} + \vec{l} \wedge \vec{m} \]
where the first term on the right, the inner product,  is a scalar (0-vector) and the second is, as we already know from above, a 2-vector. This geometric product can be extended consistently to the whole algebra to produce an associative product that encodes a variety of geometric properties of the Euclidean plane.
\item In the case of the euclidean plane, the inner product of two planes $\vec{l}_1 \cdot \vec{l}_2 := a_1 b_1 + a_2 b_2$.  Notice that the inner product is degenerate since the $0^{th}$ coordinate plays no role. This explains the notation for the algebra $\pdclal{2}{0}{1}$.
\item A \emph{k-blade} is the product of $k$ 1-vectors. 
\item We abbreviate the lowest grade part of the product of two blades $\vec{X} \vec{Y}$ as $\vec{X} \cdot \vec{Y}$. This is consistent with the above definition for the product of 1-vectors.
\item We write $\widetilde{\vec{X}}$ to denote the reversal of $\vec{X}$, the multivector obtained by reversing the order of all products of 1-vectors in the element.  This algebra involution is important in the discussion of rotations below (\Sec{sec:isom}).
\item $\eye := \e{0}\e{1}\e{2}$  represents the unit pseudoscalar, the basis of the grade-3 vectors.
%\item  For a $k$-vector $\vec{A}$ and a $l$-vector $\vec{B}$ such that $k+l=3$, we  occasionally use the fact the \[ \vec{A} \wedge \vec{B} = (\vec{A} \vee \vec{B})\eye\]  This is sometimes convenient to produce a scalar instead of a pseudoscalar. 
\end{itemize}

With this much too-brief introduction we turn now to \quot{translating} selected contents of \cite{wildberger2005} into PGA.

\subsection{Quadrance}

Given two points $A = (a_{0}, a_{1})$ and $B = (b_{0}, b_{1})$ in $\R{2}$, then the \emph{quadrance} of $A$ and $B$ is defined as the square of their euclidean distance (\cite{wildberger2005},  Ch. 5):
\[ Q(A,B) =  (a_{0} - b_{0})^{2} + (a_{1} - b_{1})^2 \]
Translating these points into the geometric algebra, we obtain  $\vec{A} = \EE{0}+a_{1}\EE{1}+a_{2}\EE{2}$ and $\vec{B} = \EE{0}+b_{1}\EE{1}+b_{2}\EE{2}$.  In this form, the quadrance can be written:
\begin{equation} \label{eqn:quadrance}
Q(\vec{A}, \vec{B}) = (\vec{A} \vee \vec{B})^{2} 
\end{equation}
 In words, the variance of $\vec{A}$ and $\vec{B}$ is the square of the  line joining the two points.

(Note: This assumes we have normalized coordinates so $a_{0} = b_{0} = 1$. Due to the degenerate nature of the euclidean metric, $\vec{A}^{2} = a_{0}^{2}$. Hence normalizing a euclidean point does not require square roots as normalizing a euclidean line does.  \Sec{sec:comparison} shows how to calculate the quadrance when the points are not necessarily normalized.)

\subsection{Spread}

Given two lines $l: a_{1}x + b_{1}y + c_{1} = 0$ and $m : a_{2}x + b_{2}y + c_{2} = 0$, then the \emph{spread} of $l$ and $m$ could be defined as the square of the sine of their angle  (\cite{wildberger2005},  Ch. 6).  The spread can be defined without reference to the \emph{sine} function, as follows.  If $A$ is the intersection of $l$ and $m$, and $B$ is some point on $A$, let $C$ be the orthogonal projection of $B$ onto $m$. (See \cite{wildberger2005}, p. 6).  Then the spread is defined as:
\[ s(l,m) := \dfrac{Q(A,B)}{Q(B,C)} = \dfrac{(a_{1}b_{2}-a_{2}b_{1})^{2}}{(a_{1}^{2}+b_{1}^{2})(a_{2}^{2}+b_{2}^{2})}\]
Translating these lines into PGA, we obtain  $\vec{l} = c_{1}\e{0}+a_{1}\e{1}+b_{1}\e{2}$ and $\vec{m} = c_{2}\e{0}+a_{2}\e{1}+b_{2}\e{2}$.  In this form, the spread can be written:
\begin{equation} \label{eqn:spread}
 s(\vec{l}, \vec{m}) = \dfrac{-(\vec{l} \wedge \vec{m})^{2}}{\vec{l}^{2}\vec{m}^{2}} 
 \end{equation}
$\vec{l} \wedge \vec{m}$ is the intersection point of $\vec{l}$ and $\vec{m}$; squaring it effectively picks off the (negative) square of the homogeneous weight of this point since $\EE{0}^{2} = -1$ and $\EE{1}^{2}=\EE{2}^{2}=0$.  The denominator of the fraction consists of the product of the squares of the two lines. If the two lines are normalized -- as they always can be -- to have norm 1, then the formula reduces to the simpler form:  $ s(\vec{l}, \vec{m}) = -(\vec{l} \wedge \vec{m})^{2}$.

\myboldhead{Cross and twist}  Notice that the concepts of cross and twist from \cite{wildberger2005} can be translated into PGA as follows:
\begin{align*} \label{eqn:spread}
 c(\vec{l}, \vec{m}) = \dfrac{(\vec{l} \cdot \vec{m})^{2}}{\vec{l}^{2}\vec{m}^{2}}  \\
 t(\vec{l}, \vec{m}) =  \dfrac{-(\vec{l} \wedge \vec{m})^{2}}{(\vec{l}\cdot\vec{m})^2}
 \end{align*}

\subsection{Comparison}
\label{sec:comparison}

Notice that the two formulae \Eq{eqn:quadrance} and \Eq{eqn:spread} exhibit a fundamental similarity which is not obvious in the coordinate-based approach.  This can be intensified if one allows non-normalized points in \Eq{eqn:quadrance}.  Then one is led to the following formula:

\begin{equation} \label{eqn:quadrance2}
Q(\vec{A}, \vec{B}) = \dfrac{(\vec{A} \vee \vec{B})^{2} }{\vec{A}^{2}\vec{B}^{2}}
\end{equation}

In this form the two formulae are, overlooking the minus sign, structurally identical. We are justified in saying that  quadrance and spread are \emph{dual} to one another.  

This duality is not obvious from the way quadrance and spread are originally defined.  Only in the context of PGA does this symmetry become clear.  Another advantage of the approach via geometric algebra  is that  the derivation of \Eq{eqn:spread} does not require the introduction of the secondary points $B$ and $C$ as in \cite{wildberger2005}. The formula is written entirely in terms of the intersection point $\vec{A} = \vec{l}\wedge \vec{m}$ and the two lines $\vec{l}$ and $\vec{m}$.

\section{Cartesian coordinate geometry}

In this section I will attempt to translate selected elements of Chapter 3 of \cite{wildberger2005} into $\pdclal{2}{0}{1}$.  The numbers in square brackets refer to the corresponding sections of \cite{wildberger2005}.

\myboldhead{Points and Lines [3.1]}
A point $\bf{P}=(x,y)$, as noted above, is mapped to the 2-vector $\EE{0} + x\EE{1} + y\EE{2}$.  In this form it is normalized and  satisfies $\bf{P}^2=-1$.  A free vector $\bf{V}$  has the form  $ x\EE{1} + y\EE{2}$. It is an \emph{ideal} point and satisfies $\bf{V}^2=0$.

A line $\langle a:b:c \rangle$ is mapped to the 1-vector $l := c\e{0} + a\e{1} +b\e{2}$. We can always normalize this line to satisfy $l^2 = -1$ whenever $a^2+b^2\neq0$.  In the latter case we say the line is \emph{ideal}, otherwise we say the line is \emph{proper}, see \cite{gunn2011}.  \cite{wildberger2005} refers to \emph{null}, resp. \emph{non-null} lines.  A point $\vec{P}$ lies on a line $\vec{m}$ precisely when $\vec{P} \wedge \vec{m} = \vec{P} \vee \vec{m} = 0$.

\myboldhead{Collinearity and concurrence [3.2]} The line $\vec{l}$ through the two points $\vec{A}_1$ and $\vec{A}_2$ is given by $\vec{l} = \vec{A}_1 \vee \vec{A}_2$.  Three points are collinear precisely when $\vec{A}_1 \vee \vec{A}_2 \vee \vec{A}_3  = 0$.  Similarly, three lines are concurrent if and only if $\vec{l}_1 \wedge \vec{l}_2 \wedge \vec{l}_3  = 0$. 

\myboldhead{Parallel and perpendicular lines [3.3]} Two distinct lines are parallel when $\vec{l}_1 \wedge \vec{l}_2 \wedge \e{0} = 0$, that is, the intersection point lies on the ideal line. This is equivalent to the condition $(\vec{l}_1 \wedge \vec{l}_2 )^2 = 0$, that is, the spread of the two lines vanishes. Two lines are perpendicular exactly when $\vec{l}_1 \cdot \vec{l}_2 = 0$.  

\myboldhead{Altitude and parallel to a line [3.4]}. Given a point $\vec{A}$ and a line $\vec{l}$, $\vec{l}_\perp := \vec{A}\cdot\vec{l}$ represents the grade-1 part of the product $\vec{A} \vec{l}$: it is the line through $\vec{A}$ perpendicular to $\vec{l}$, or the \emph{altitude} from $\vec{A}$ to $\vec{l}$.   Multiplying this by $\vec{A}$ produces the \emph{parallel} through $\vec{A}$ to $\vec{l}$. $\vec{l}_\parallel := \vec{A}(\vec{A}\cdot\vec{l})$.  With the same notation as above, the \emph{foot} of the altitude  is given by the intersection of $\vec{l}$ with $\vec{l}_\perp$:  $\vec{l} \wedge \vec{l}_\perp$.

\myboldhead{Perpendicular bisector [3.8]}  Given two  (normalized) points $\vec{A}_1$ and $\vec{A}_2$, the perpendicular bisector of the line $\vec{A}_1 \vec{A}_2$ is given by $(\vec{A}_1 \vee \vec{A}_2)\cdot(\vec{A}_1+ \vec{A}_2)$.  Here the first term is the line joining the points, while the second is the midpoint of their segment.  In this context, the inner product $\cdot$ of a line and a point was defined in \ref{sec:crashcourse} as the lowest-grade part of the geometric product of the line and the point.  In this case it is the line through the point perpendicular to the line (for details see \cite{gunnFull2010}).
%It is remarkable how the themes of this section map onto the structure of the algebra.

\subsection{Triangles [3.5]}

Here we deal with the remarks in Chapter 3 dealing with triangles. We assume that all given points are normalized to have homogeneous coordinate 1. Given a triangle with vertices  $\vec{A}_1$, $\vec{A}_2$, and $\vec{A}_3$.  Then the sides can be defined $\vec{a}_i := \vec{A}_j \vee \vec{A}_k$.  .  Let $\tarea$ represent the signed area of the triangle $\Delta \vec{A}_1\vec{A}_2\vec{A}_3$.   Then:
\begin{align} \label{eqn:triarea}
\vec{A}_1 \vee \vec{A}_2 \vee \vec{A}_3 &=  2 \tarea 
\end{align}
The various products involving the sides can also be calculated:
\begin{align} \label{eqn:asqr}
\vec{a}_i^2  &= Q(A_j, A_k) \\ \label{eqn:aiwedgeaj}
\vec{a}_i \wedge \vec{a}_j &= 2 \tarea \vec{A}_k
\end{align}
In words: the wedge of two lines of the triangle is the included vertex times twice the area of the triangle.  We can now derive the formula for the triple wedge of the sides:
\begin{align}
\vec{a}_1 \wedge \vec{a}_2 \wedge \vec{a}_3 &= (\vec{a}_1 \wedge \vec{a}_2) \wedge \vec{a}_3   \\
&=  2 \tarea \vec{A}_3  \wedge \vec{a}_3 \\
&= 2 \tarea (\vec{A}_3  \vee \vec{a}_3)\eye \\
&= 2 \tarea (\vec{A}_1  \vee \vec{A}_2 \vee \vec{A}_3)\eye \\
&=  4 \tarea^2 \eye
\end{align}

In the above derivation, we have used associativity, \Eq{eqn:aiwedgeaj},  duality of $\vee$ and $\wedge$, definition of $\vec{a}_i$, and \Eq{eqn:triarea}.

\myboldhead{Centroid}  Given a triangle $\vec{A}_i$ for $i\in\{1,2,3\}$. Then the median $\vec{m}_i$ is defined as the joining line of $\vec{A}_i$ with the midpoint $\vec{M}_i$ of the opposite side.  A little reflection shows that this is given by $\vec{m}_i = (\vec{A}_{j}+\vec{A}_{k}) \vee \vec{A}_i$ where $(i,j,k)$ is a cyclic permutation of $(1,2,3)$.  To show that the three medians are concurrent, we calculate the intersection of a pair of medians $\vec{m}_i  \wedge \vec{m}_j$ and show it is independent of $i$ and $j$.  

\begin{align*}
\vec{m}_i \wedge \vec{m}_j  &= (\vec{A}_j \vee \vec{A}_i +\vec{A}_k \vee \vec{A}_i) \wedge (\vec{A}_k \vee \vec{A}_j +\vec{A}_i \vee \vec{A}_j) \\
&= (-\vec{a}_k + \vec{a}_j) \wedge (-\vec{a}_i + \vec{a}_k) \\
&= \vec{a}_k \wedge \vec{a}_i+ \vec{a}_i\wedge \vec{a}_j + \vec{a}_j \wedge \vec{a}_k
\end{align*}

 For more on handling triangle centers within $\pdclal{2}{0}{1}$, including a treatment of the Euler line, see \cite{gunnFull2010}. p. 17.

\myboldhead{Thales theorem} There is a \quot{second} theorem attributed to Thales:

\emph{ If a straight line is drawn parallel to one of the sides of a triangle $\bf{A_1 A_2 A_3}$, then it cuts the other sides of the triangle, or these produced, proportionally; and, if the sides of the triangle, or the sides produced, are cut proportionally, then the line joining the points of section is parallel to the remaining side of the triangle.
}

\textbf{Proof}:  The first half of the theorem can be easily proved using similar triangles (also the converse?).  For the converse: the condition that the line cuts the two sides $\bf{A_1 A_2}$ and $A_1 A_3$ proportionally is that there is a real $\lambda \neq 0$ such that the respective intersection points $\bf{B_2}$ and $\bf{B_3}$ satisfy: 
\begin{eqnarray*}
\vec{B}_2 := (1-\lambda)\vec{A}_1 + \lambda \vec{A}_2 \\
\vec{B}_3 := (1-\lambda)\vec{A}_1 + \lambda \vec{A}_3
\end{eqnarray*}
We show that the lines $\vec{l}_A := \vec{A}_2 \vee \vec{A}_3$ and $\vec{l}_B := \vec{B}_2 \vee \vec{B}_3$ are parallel by showing that their intersection point is ideal:
 \begin{align*}
(\vec{l}_A \wedge \vec{l}_B) &=  ( \vec{A}_2 \vee \vec{A}_3) \wedge \vec{B}_2 \vee \vec{B}_3 \\
&= ( \vec{A}_2 \vee \vec{A}_3) \wedge ( (1-\lambda)\vec{A}_1 + \lambda \vec{A}_2) \vee ( (1-\lambda)\vec{A}_1 + \lambda \vec{A}_3)  \\
&= ( \vec{A}_2 \vee \vec{A}_3) \wedge ((1-\lambda)\lambda ( \vec{A}_1\vee  \vec{A}_3 + \vec{A}_2 \vee \vec{A}_1) + \lambda^2 \vec{A}_2 \vee \vec{A}_3) \\
&= (\vec{A}_2 - \vec{A}_3)(1-\lambda)\lambda 2 \tarea
\end{align*}

The final expression is an ideal point since the difference of two points is ideal (coefficient of $\EE{0}$ vanishes). QED.

 \section{ Reflections  and other isometries}
 \label{sec:isom}
 Here we handle the remarks in Chapter 4, mostly concerned with euclidean isometries.
 
% \subsection{Reflections}
 One of the great advantages of geometric algebra is the uniform representation of isometries.  The reflection in the line $\vec{m}$ is represented by the \emph{sandwich operation} $\vec{X} \rightarrow \vec{m} \vec{X} \vec{m}$.  Here the grade of $\vec{X}$ is arbitrary: the same expression works for 0-, 1-, 2-, and 3-vectors $\vec{X}$.  
 
Reflections, as is well known, generate the complete euclidean group $\Euc{2}$.  In geometric algebra the product of two reflections can be written as a \quot{thick}  sandwich:  \[\vec{m}(\vec{l} \vec{X} \vec{l})\vec{m} = (\vec{m} \vec{l})\vec{X} (\vec{l}\vec{m})\]  where we have used associativity to move the parentheses. This sandwich represents a rotation around the point $\vec{P} = \vec{m} \wedge \vec{l}$.  (If the two lines are parallel, $\vec{P}$ is ideal and the isometry is a translation.) 
Which rotation is it?

If $\vec{l}$ and $\vec{m}$ are normalized to $\vec{m}^2=\vec{l}^2 = 1$, then \[ \vec{R} := \vec{m} \vec{l} = \cos{\theta} + \sin{\theta}\vec{P}\] where $\theta$ is the angle between the lines $\vec{l}$ and $\vec{m}$.  The rotation then takes the form $\vec{R} \vec{X} \widetilde{\vec{R}}$. In general one has  \[ \vec{R}  = \sqrt{ \vec{l}^{2}\vec{m}^{2}} (\cos{\theta} + \sin{\theta}\vec{P})\] The resulting rotation is through the angle $2\theta$.  An analogous formula involving translation distance holds for translations.  Note that  $\vec{P}$ is also normalized to have homogeneous coordinate 1 in this expression.  

This analysis makes clear why it is useful to work with normalized lines and points when working with isometries; on the other hand, the discussion of triangles above shows that on other occasions one does not want to normalize.

 \section{Conclusion}
 
 This translation of a small sampling from Part I: \emph{Preliminaries} of \cite{wildberger2005} demonstrates that the algebra $\pdclal{2}{0}{1}$ is a natural environment for presenting these results.  Statements and proofs are often shorter and more direct than the analogous calculations involving coordinates.  One also notices otherwise hidden connections, such as the underlying, exact duality of quadrance and spread. I am not aware of any result in universal geometry which cannot be translated to PGA; the paucity of results presented here reflects a lack of time rather than any intrinsic obstacles to translation.  
 
 This offers promising outlook for further work, either in the form of extending the results to  Part II of  \cite{wildberger2005}, but also for other work of the same author devoted to non-euclidean planar geometry.  Here the algebras $\pdclal{3}{0}{0}$ and $\pdclal{2}{1}{0}$ promise to provide the same faithful representation for elliptic (spherical) and hyperbolic geometry, resp.  The corresponding reference is then the thesis \gTh, Chapter 6.
 
% \begin{thebibliography}{Gun11b}
%
%\bibitem[Gun11a]{gunnThesis}
%Charles Gunn.
%\newblock {\em Geometry, Kinematics, and Rigid Body Mechanics in Cayley-Klein
%  Geometries}.
%\newblock PhD thesis, Technical University Berlin, 2011.
%\newblock \url{http://opus.kobv.de/tuberlin/volltexte/2011/3322}.
%
%\bibitem[Gun11b]{gunn2011}
%Charles Gunn.
%\newblock On the homogeneous model of euclidean geometry.
%\newblock In Leo Dorst and Joan Lasenby, editors, {\em A Guide to Geometric
%  Algebra in Practice}, chapter~15, pages 297--327. Springer, 2011.
%\newblock
%  \url{http://page.math.tu-berlin.de/~gunn/Documents/Papers/GAInPracticeCh15Gunn.pdf}.
%
%\bibitem[Gun11c]{gunnFull2010}
%Charles Gunn.
%\newblock On the homogeneous model of euclidean geometry: Extended version.
%\newblock {\em \url{http://arxiv.org/abs/1101.4542}}, 2011.
%
%\bibitem[Gun17a]{gunn2017b}
%Charles Gunn.
%\newblock Doing euclidean plane geometry using projective geometric algebra.
%\newblock {\em Advances in Applied Clifford Algebras}, 27(2):1203--1232, 2017.
%
%\bibitem[Wil05]{wildberger2005}
%N.~J. Wildberger.
%\newblock {\em Divine Proportions}.
%\newblock Wild Egg, 2005.
%

%\end{thebibliography}

\bibliography{GunnRef}
\bibliographystyle{alpha}

\end{document}

%% file: GunnMac.tex
\usepackage{xspace}
\usepackage{helvet}
\usepackage{courier}
\usepackage{fancyhdr}
\usepackage{type1cm}         
\usepackage{subeqnarray}
\usepackage{graphicx}        % standard LaTeX graphics tool
\usepackage{enumerate}
\usepackage{amsmath, latexsym}
\usepackage{array}
\usepackage{tabularx}
\usepackage{url}
\usepackage{relsize}
\usepackage{ifthen}
\usepackage{xcolor}
\usepackage[newitem,newenum]{paralist}
 
\def\isBeamer{false}

\newcommand{\mychapterskip}{\ifthenelse{\equal{\isMPK}{true}}
{\vspace{-.8in}}{}}

\newcommand{\fullversioncolor}{\color{black!50!blue}}
\newcommand{\versions}[2]{\ifthenelse{\equal{\isfullversion}{true}}{{\fullversioncolor#1}}{#2}}
\renewcommand{\vec}[1]{\mathbf{#1}}
\newcommand{\quot}[1]{``#1''}

\newcommand{\R}[1]{\ensuremath{\mathbb{R}^{#1}\xspace}}
\newcommand{\Euc}[1]{\ensuremath{\mathbf{E}^{#1}\xspace}}

\newcommand\proj[1]{\ensuremath{\mathbf{P}(#1)}\xspace}

\newcommand{\e}[1]{\vec{e}_{#1}}
\newcommand{\EE}[1]{\vec{E}_{#1}}

\newcommand{\eye}{\vec{I}}

\newcolumntype{Y}{X}

\newcommand{\pdclal}[3]{\proj{\mathbb{R}^*_{#1,#2,#3}}\xspace}

\newcommand{\myboldhead}[1]{\vspace{0in}\hspace{-.15in}\textbf{#1.}}

\newcommand{\app}[1]{\Sec{sec:J}}

\definecolor{mypurple}{rgb}{1,0,1}

\definecolor{mygreen}{rgb}{0, .5, 0}

\newcommand{\Eq}[1]{(\ref{#1})}

\newcommand{\Sec}[1]{Sect.~\ref{#1}}

\ifthenelse{\equal{\isBeamer}{true}}
{

}
{

}

\newcommand{\mycorrection}[1]{}

%% file: rationalTrig.bbl
\begin{thebibliography}{Gun17b}

\bibitem[Gun11a]{gunnThesis}
Charles Gunn.
\newblock {\em Geometry, Kinematics, and Rigid Body Mechanics in Cayley-Klein
  Geometries}.
\newblock PhD thesis, Technical University Berlin, 2011.
\newblock \url{http://opus.kobv.de/tuberlin/volltexte/2011/3322}.

\bibitem[Gun11b]{gunn2011}
Charles Gunn.
\newblock On the homogeneous model of euclidean geometry.
\newblock In Leo Dorst and Joan Lasenby, editors, {\em A Guide to Geometric
  Algebra in Practice}, chapter~15, pages 297--327. Springer, 2011.

\bibitem[Gun11c]{gunnFull2010}
Charles Gunn.
\newblock On the homogeneous model of euclidean geometry: Extended version.
\newblock {\em \url{http://arxiv.org/abs/1101.4542}}, 2011.

\bibitem[Gun17a]{gunn2017b}
Charles Gunn.
\newblock Doing euclidean plane geometry using projective geometric algebra.
\newblock {\em Advances in Applied Clifford Algebras}, 27(2):1203--1232, 2017.
\newblock \url{http://arxiv.org/abs/1501.06511 }.

\bibitem[Gun17b]{gunn2017a}
Charles Gunn.
\newblock Geometric algebras for euclidean geometry.
\newblock {\em Advances in Applied Clifford Algebras}, 27(1):185--208, 2017.
\newblock \url{http://arxiv.org/abs/1411.6502 }.

\bibitem[Wil05]{wildberger2005}
N.~J. Wildberger.
\newblock {\em Divine Proportions}.
\newblock Wild Egg, 2005.

\end{thebibliography}
